\documentclass[12pt]{article}

\usepackage{tikz}
\usepackage{subfigure}
\usepackage[english]{babel}
\usepackage{mathrsfs}
\usepackage[center]{caption2}
\usepackage{amsfonts,amssymb,amsmath,latexsym}
\usepackage{multirow}
\usepackage{amsthm}

\def\g_p{\gamma_{pr}}

\newtheorem{thm}{Theorem}

\newtheorem{lem}[thm]{Lemma}

\newtheorem{quest}[thm]{Question}
\newtheorem{claim}{Claim}
\newtheorem{fact}{Fact}

\begin{document}
\title{$H$-Decomposition of $r$-graphs when $H$ is an $r$-graph with exactly $k$ independent edges
\thanks{The work was supported by NNSF of China (No. 11671376) and  NSF of Anhui Province (No. 1708085MA18).}
}
\author{Xinmin Hou$^a$,\quad Boyuan Liu$^b$, \quad Hongliang Lu$^c$\\
\small $^{a,b}$ Key Laboratory of Wu Wen-Tsun Mathematics\\
\small School of Mathematical Sciences\\
\small University of Science and Technology of China\\
\small Hefei, Anhui 230026, China.\\
\small $^c$ School of Mathematics and Statistics\\
\small Xi'an Jiaotong University\\
\small Xi'an, Shaanxi 710049, China\\
}

\date{}


\maketitle

\begin{abstract}
Let $\phi_H^r(n)$ be the smallest integer such that, for all $r$-graphs $G$ on $n$ vertices, the edge set $E(G)$ can be partitioned into at most $\phi_H^r(n)$ parts, of which every part either is a single edge or forms an $r$-graph isomorphic to $H$.
The function $\phi^2_H(n)$ has been well studied in literature, but for the case $r\ge 3$, the problem that determining the value of $\phi_H^r(n)$ is widely open. Sousa (2010)
gave an asymptotic value of $\phi_H^r(n)$ when $H$ is an $r$-graph with exactly 2 edges, and determined the exact value of $\phi_H^r(n)$ in some special cases.
In this paper,  we first give the exact value of $\phi_H^r(n)$ when $H$ is an $r$-graph with exactly 2 edges, which improves Sousa's result. Second we determine the exact value of $\phi_H^r(n)$ when $H$ is an $r$-graph consisting  of exactly $k$ independent edges.
\end{abstract}

\section{Introduction}
Given two $r$-graphs $G$ and $H$, an {\em $H$-decomposition} of $G$ is a partition of the edge set of $G$  such that each part is either a single edge or forms an $r$-graph isomorphic to $H$. The minimum number of parts in an $H$-decomposition of $G$ is denoted by $\phi_H^r(G)$. The $H$-decomposition number $\phi_H^r(n)$ is defined as
\begin{equation*}
  \phi_H^r(n)=\max\{\phi_H^r(G) : G \mbox{ is an $r$-graph with } |V(G)|=n\}.
\end{equation*}
{ An $r$-graph $G$ with $\phi_H^r(G)=\phi_H^r(n)$ is called an {\it extremal graph} of $H$.}

For the case $r=2$, we omit the index $2$ for short in the paper, for example we write  graph for  2-graph, and $\phi_H(n)$ for $\phi_H^2(n)$. The function $\phi_H(n)$ has been well studied in literature by many researchers. The first exact value of $\phi_H(n)$ when $H=K_3$ was given by Erd\H{o}s, Goodman and P\'{o}sa~\cite{intersections} in 1966,  where $K_k$ is the complete graph on $k$ vertices.  Ten years later, Bollob\'{a}s~\cite{Bollobas} generalized the result to $H=K_k$, $k\ge 3$. Much more exact values of $\phi_H(n)$ can be found in the survey by Sousa~\cite{Survy} in 2015. Recently, Hou, Qiu and Liu determined the exact values of $\phi_H(n)$ when $H$ is a graph consisting of $k$ complete graphs of order at least 3 which intersect in exactly one common vertex~\cite{Dec for Fkr} and when $H$ is a graph consisting of $k$ cycles of odd lengths which intersect in exactly one common vertex~\cite{Dec Hst}.
An asymptotic value of the function $\phi_H(n)$ was given by Pikhurko and Sousa~\cite{P&S}, and lately the value was improved  by Allen, B\"{o}ttcher, and Person~\cite{ImprovedError}.

For the case $r\ge 3$, the study of the function $\phi_H^r(n)$ is widely open.  Sousa~\cite{origin}
gave an asymptotic value of $\phi_H^r(n)$ when $H$ is an $r$-graph consisting of 2 edges, and determined the exact value of $\phi_H^r(n)$ in the special cases where the two edges of $H$ intersect exactly 1, 2 and $r-1$ vertices.
In this paper, we first generalize Sousa's result in~\cite{origin}, that is we obtain the exact value of $\phi_H^r(n)$ when $H$ is an $r$-graph consisting of exactly 2 edges. Second we focus on the case that $H$ is the $r$-graph consisting  of exactly $k$ independent edges, and we determine the exact value of $\phi_H^r(n)$ in this case.

Given positive integer $n$, $r$ and $k$ with $n\ge r\ge 2$, let $K_n^r$ be a complete $r$-graph on $n$ vertices and let ${K_n^r-\ell e}$ be a graph obtained from $K_n^r$ by deleting $\ell$ edges from $K_n^r$, where $\ell$ is an  integer so that  $0\le \ell\le k-1$ and $e(K_n^r-\ell e)\equiv k-1\pmod k$.
Note that $\ell$ is determined uniquely by $n, r$ and $k$. Let $\mathcal{K}_{n,k}^r$ be the family of $r$-graphs $K_n^r-\ell e$. The followings are our main results.
\begin{thm}\label{main1}
{ Given integers $r,k$ satisfying $0\le k\leq r-1$. Let $H$ be an $r$-graph consisting of exactly 2 edges which intersect k vertices.}
If {$n\geq2r-k$}, then $\phi_H^r(n)=\left\lceil\frac{1}{2}\binom{n}{r}\right\rceil$. Moreover,
graphs $G\in\mathcal{K}_{n,2}^r$  and $G=K_n^r$ if ${n\choose r}\equiv 0\pmod2$ are extremal graphs of $H$.


\end{thm}

Theorem~\ref{main1} improves Sousa's result in~\cite{origin}.

\begin{thm}\label{main2}
Given integers $k\ge 1$, $r\ge 2$ and {$n_0=kr(k+r-2)+2r-1$}, let $H$ be an $r$-graph on $n$ vertices consisting of exactly $k$ independent edges. If $n\geq n_0$ then
\begin{equation*}
\phi_H^r(n)=\left\{\begin{array}{ll}
                     \left\lfloor \frac{1}{k}\binom{n}{r}\right\rfloor+k-1, & \mbox{if } {n\choose r}\equiv k-1\pmod k;\\
                     \left\lfloor \frac{1}{k}\binom{n}{r}\right\rfloor+k-2, & \mbox{otherwise}.
                   \end{array}
           \right.
\end{equation*}
Furthermore, $G$ is an extremal graph of $H$ if and only if $G\in\mathcal{K}_{n,k}^r$ or $G=K_n^r$ if $\binom{n}{r}\equiv k-2\pmod k$.


\end{thm}

The proofs of Theorems~\ref{main1} and~\ref{main2} will be given in Sections 2 and 3, respectively. Before giving the proofs, we first introduce some definitions and {notation}. Let $H$ be an $r$-graph with vertex set $V(H)$ and $E(H)$. For a vertex $v\in V(H)$, the {\it degree} of $v$, denoted by $d_H(v)$,  is the number of edges of $H$ containing $v$, and the {\it minimum degree} of $H$ is denoted by $\delta(H)$.
 The {\it matching number} of $H$ is the maximum number of independent edges in $H$. We write $e(H)$ for the number of edges of $H$, that is $e(H)=|E(H)|$.

\section{Proof of Theorem~\ref{main1}}
We need some basic facts in algebraic graph theory.
A graph $G$ is called {\it vertex-transitive} if its automorphism group acts transitively on $V(G)$. Given nonnegative integers $n, r$ and $k$, let $J(n,r,k)$ be the graph  with vertex set $E(K_n^r)$, where two vertices are adjacent if and only if their intersection has size $k$. For $n\ge r$, the graphs $J(n,r,r-1)$ and $J(n,r,0)$ are known as the Johnson graphs and the Kneser graphs, respectively.

\begin{fact}[See page 9 and page 35 in~\cite{vertex transitive}]\label{FACT: Johnson-graph}
(1) $J(n,r,k)$ has $n\choose r$ vertices, and each vertex has degree ${r\choose k}{n-r\choose r-k}$.

(2)  The graphs $J(n,r,k)$ are vertex-transitive.

(3) If $n\ge r\ge k$, $J(n,r,k)\cong J(n,n-r,n-2r+k)$.

\end{fact}

\begin{lem}[Theorem 3.5.1 in~\cite{vertex transitive}]\label{vertex transitive}
If $G$ is a connected vertex-transitive graph, then $G$ has a matching that misses at most one vertex.
\end{lem}


\begin{lem}[Theorem 2.3 in~\cite{origin}]\label{LEM: upper-bound}
 Let $H$ be a fixed $r$-graph with 2 edges and $G$ an $r$-graph with $n$ vertices. Then $\phi_H^r(G)\le\phi^r_H(K^r_n)$.
\end{lem}

\noindent{\bf Proof of Theorem~\ref{main1}:} Lemma~\ref{LEM: upper-bound} implies that $\phi_H^r(n)=\phi_H^r(K_n^r)$. So, to prove the result, it is sufficient to show that $\phi_H^r(K_n^r)=\lceil\frac{1}{2}\binom{n}{r}\rceil$. Clearly, $\phi_H^r(K_n^r)\ge \lceil\frac{1}{2}\binom{n}{r}\rceil$ as $e(H)=2$.
To prove $\phi_H^r(K_n^r)\le \lceil\frac{1}{2}\binom{n}{r}\rceil$, it is sufficient to find an $H$-decomposition of $K_n^r$ with $\lceil\frac{1}{2}\binom{n}{r}\rceil$ parts.

By the definition of $J(n,r,k)$, $K_n^r$ has an $H$-decomposition with $\lceil\frac{1}{2}\binom{n}{r}\rceil$ parts is equivalent to
the statement that $J(n,r,k)$ has a matching missing at most one vertex. By (3) of Fact~\ref{FACT: Johnson-graph}, we may assume $n\ge 2r$.
If $k=0$ and $n=2r$, then $J(2r,r,0)$ consists of exactly $\frac{1}{2}\binom{n}{r}$ independent edges, so the statement holds.

Now we assume $k>0$ or $n>2r$.    By (2) of Fact~\ref{FACT: Johnson-graph}, $J(n,r,k)$ is vertex-transitive. Hence, by Lemma~\ref{vertex transitive}, to show $J(n,r,k)$ has a matching missing at most one vertex, it is sufficient to show that $J(n,r,k)$ is connected.
{That is, we need to show that {any pair of vertices $e, f$ of $J(n,r,k)$ are connected}. Suppose $|e\cap f|=i$. We prove $e$ and $f$ are connected by induction on $i$. If $i=r$, then statement is trivial. If $i=r-1$,  assume $e=\{1,2,\ldots,r-1,r\}$ and $f=\{1,2,\ldots,r-1,r+1\}$. Then $ehf$ with $h=\{1,2,\ldots,k,r+2,\ldots,2r+1-k\}$ is a walk connecting $e$ and $f$ in $J(n,r,k)$. So the result is true for the base case. Now suppose $i<r-1$ and the statement is true for any large $i$. By symmetry, one could assume that $e=\{1,2,\ldots,r-1,r\}$ and $f=\{1,2,\ldots,i,r+1,\ldots, 2r-i\}$. Let $h=\{1,\ldots,i,i+1,r+1,~\ldots,2r-i-1\}$. Then $|e\cap h|=i+1>i$ and $|f\cap h|=r-1>i$. By induction hypothesis, $e$ and $h$ (resp. $f$ and $h$) are connected in $J(n,r,k)$. By the transitivity of connectivity, $e$ and $f$ are connected in $J(n,r,k)$.
\medskip\rule{1mm}{2mm}



\section{Proof of Theorem \ref{main2}}
We need two known theorems to prove our result.
Given graphs $G$ and $H$, we say $G$ has an {\it $H$-factor} if $G$ contains $\left\lfloor \frac{|V(G)|}{|V(H)|}\right\rfloor$ vertex-disjoint copies of $H$.

\begin{thm}[Hajnal, Szemer\'edi~\cite{k-factor}]\label{factor}\label{k-factor}
Let $k$ be a positive integer. If $G$ is a graph on $n$ vertices with minimum degree
\begin{equation*}
 \delta(G)\geq(1-\frac{1}{k})n,
\end{equation*}
then $G$ has a $K_k$-factor.
\end{thm}

\begin{thm}[Frankl \cite{matching}]\label{matching}
If $H$ is an $r$-graph on $n$ vertices with matching number of size $k$ and $n\geq (2k+1)r-k$, then
\begin{equation*}
  e(H)\leq \binom{n}{r}-\binom{n-k}{r}.
\end{equation*}
\end{thm}

\noindent{\bf Proof of Theorem~\ref{main2}:}
Let $G$ be an $r$-graph on $n\geq n_0$ vertices with $\phi_H(G)=\phi_H(n)$. Let $p_H(G)$ denote the maximum number of pairwise edge-disjoint copies of $H$ in $G$. Then we have
\begin{equation}\label{EQN: e1}
\phi_H(G)=e(G)-(k-1)p_H(G).
\end{equation}
If we remove the edges of $p_H(G)$ pairwise edge-disjoint copies of $H$ from $G$, then we obtain an $H$-free graph, that is a graph with matching number at most $k-1$. Hence by Theorem \ref{matching}, we have
\begin{equation}\label{EQN: e2}
\binom{n}{r}-\binom{n-k+1}{r}\geq e(G)-k p_H(G).
\end{equation}
On the other hand,
\begin{equation}\label{EQN: e3}
\phi_H(G)\geq\phi_H(K_n^r)\geq \frac 1k{\binom{n}{r}}.
\end{equation}
From   (\ref{EQN: e1}), (\ref{EQN: e2}) and (\ref{EQN: e3}), we have
\begin{equation}\label{EQN: e4}
e(G)\geq \binom{n}{r}-(k-1)\left[\binom{n}{r}-\binom{n-k+1}{r}\right].
\end{equation}

Now we define an auxiliary graph $L_G$ as follows: let $V(L_G)=E(G)$ and two vertices $e_1, e_2\in V(L_G)$ is adjacent if and only if $e_1\cap e_2=\emptyset$ in $G$. Hence the edge set of a copy of $H$ in $E(G)$ induces a clique of order $k$ in $L_G$. Therefore, a collection of edge-disjoint copies of $H$ in $G$ corresponds to a collection of vertex-disjoint $K_k$ in $L_G$.

\begin{claim}\label{claim}
$L_G$ has a $K_k$-factor. In particular, $L_{K_n^r-\ell e}$ has a $K_k$-factor for every $\ell$ with $0\le\ell\le k-1$ .
\end{claim}

\noindent{\it Proof of the claim:} By definition of $L_G$, we have
\begin{equation*}
\delta(L_G)\geq \binom{n-r}{r}-\left[\binom{n}{r}-e(G)\right].
\end{equation*}
By Theorem~\ref{factor}, it suffices to show that
\begin{equation*}
\binom{n-r}{r}-\left[\binom{n}{r}-e(G)\right]\geq(1-\frac{1}{k})e(G),
\end{equation*}
that is
\begin{equation}\label{EQN: e5}
e(G)\geq k\left[\binom{n}{r}-\binom{n-r}{r}\right].
\end{equation}
To show (\ref{EQN: e5}), by (\ref{EQN: e4}) it suffices to show
\begin{equation*}
\binom{n}{r}-(k-1)\left[\binom{n}{r}-\binom{n-k+1}{r}\right]\geq k\left[\binom{n}{r}-\binom{n-r}{r}\right],
\end{equation*}
{that is, we need to} show
\begin{equation}\label{EQN: e6}
k\binom{n-r}{r}+(k-1)\binom{n-k+1}{r}\geq (2k-2)\binom{n}{r}.
\end{equation}
By the inequality
\begin{equation*}
\frac{\binom{n-t}{r}}{\binom{n}{r}}\geq \left(\frac{n-t-r+1}{n-r+1}\right)^r\geq 1-\frac{r t}{n-r+1}\,\,,  ({r\ge t}\ge 0)
\end{equation*}
and $n\geq n_0=kr(k+r-2)+2r-1$,  it can be easily check that  (\ref{EQN: e6}) holds. This completes the proof of the claim.

\vspace{5pt}
Now suppose $e(G)\equiv i\pmod k$ and $e(G)=t k+i\, (0\leq i\leq k-1)$ for some $t\le \lfloor\frac 1k{n\choose r}\rfloor$. By Claim~\ref{claim}, $p_H(G)=t$ and hence by (\ref{EQN: e1}), we have   $\phi_H(G)=t+i$. In particular,
\begin{equation*}
\phi_H^r(K_n^r-\ell e)=\left\{\begin{array}{ll}
                     \left\lfloor \frac{1}{k}\binom{n}{r}\right\rfloor+k-1, & \mbox{if } {n\choose r}\equiv k-1\pmod k\\
                     \left\lfloor \frac{1}{k}\binom{n}{r}\right\rfloor+k-2, & \mbox{otherwise}
                   \end{array}
           \right..
\end{equation*}

{If ${n\choose r}\equiv k-1\pmod k$, then $\phi_H^r(n)=\phi_H^r(G)=t+i\le \left\lfloor \frac 1k{\binom{n}{r}}\right\rfloor+k-1$, and the equality holds if and only if $G=K_n^r \in\mathcal{K}_{n,k}^r$. Otherwise, $\phi_H^r(n)=\phi_H^r(G)=t+i\le \left\lfloor \frac 1k{\binom{n}{r}}\right\rfloor+k-2$, the equality holds if and only if $t=\lfloor \frac 1k{\binom{n}{r}}\rfloor-1$ and $i=k-1$ or $t=\lfloor \frac 1k{\binom{n}{r}}\rfloor$ and $i=k-2$, in the former case $G\in\mathcal{K}_{n,k}^r$ and in the latter case it  happens if and only if $G=K_n^r$ and ${n\choose r}\equiv k-2\pmod k$.}\quad \rule{1mm}{2mm}


\section{Concluding Remarks}

In this paper we determine the exact value of of $\phi_H^r(n)$ when $H$ is an $r$-graph consisting of exactly 2 edges or  consisting of  exactly $k$ INDEPENDENT edges. { We believe that Theorem~\ref{main2} still holds when $H$ consists of exactly $k$ edges which intersect the same set of size $i$ $(0\leq i\leq r-1)$, we leave this as an open problem.
\begin{quest}
{Is the following statement true?} Given integers $k\ge 1$, $r\ge 2$, let $H$ be an $r$-graph consisting of exactly $k$ edges which intersect the same set of size $i$ ($0\leq i\leq r-1$). If $n$ is sufficiently large, then
\begin{equation*}
\phi_H^r(n)=\left\{\begin{array}{ll}
                     \left\lfloor \frac{1}{k}\binom{n}{r}\right\rfloor+k-1, & \mbox{if } {n\choose r}\equiv k-1\pmod k;\\
                     \left\lfloor \frac{1}{k}\binom{n}{r}\right\rfloor+k-2, & \mbox{otherwise}.
                   \end{array}
           \right.
\end{equation*}
\end{quest}
}

\end{document}